\documentclass[10pt,a4paper]{article}
\usepackage{amsmath}
\usepackage{amssymb}
\usepackage{amsthm}
\usepackage{amsfonts}
\usepackage{enumerate}
\usepackage{pstricks}
\usepackage{pst-plot}
\usepackage{pst-poly}
\usepackage{graphics} 
\usepackage{graphicx}
\usepackage{pgfpages}
\usepackage{url}
\usepackage{booktabs}
\usepackage{verbatim}
\usepackage[margin=1.3in]{geometry}
\newcommand{\tluste}[1]{\mbox{\mathversion{bold}$ #1 $}}

\newcommand{\mace}[1]{{{#1}}}

\newcommand{\omace}[1]{\mbox{$\overline{\mace{#1}}$}} 
\newcommand{\umace}[1]{\mbox{$\underline{\mace{#1}}$}} 
\newcommand{\imace}[1]{\mbox{$\tluste{#1}$}}

\def\Mid#1{{#1^c}}
\def\Rad#1{{#1^\Delta}}

\newcommand{\onum}[1]{\mbox{$\overline{{#1}}$}} 
\newcommand{\unum}[1]{\mbox{$\underline{{#1}}$}}

\newcommand{\ivr}[1]{\mbox{$\tluste{#1}$}} 

\newcommand{\inum}[1]{\mbox{$\tluste{#1}$}}

\newcommand{\R}[0]{{\mathbb{R}}}

\newcommand{\mmid}[0]{;\,}		

\DeclareMathOperator{\diag}{diag}	

\def\nref#1{$(\ref{#1})$}

\newtheorem{proposition}{Proposition}

\theoremstyle{definition}

\newtheorem{example}{Example}

\begin{document}

\title{The Effect of Hessian Evaluations in the Global Optimization $\alpha$BB Method}

\author{
  Milan Hlad\'{i}k\footnote{
Charles University, Faculty  of  Mathematics  and  Physics,
Department of Applied Mathematics, 
Malostransk\'e n\'am.~25, 11800, Prague, Czech Republic, 
e-mail: \texttt{milan.hladik@matfyz.cz}
}
}

\date{\today}
\maketitle

\begin{abstract}
We consider convex underestimators that are used in the global optimization $\alpha$BB method and its variants. The method is based by augmenting the original nonconvex function by a relaxation term that is derived from an interval enclosure of the Hessian matrix. In this paper, we discuss the advantages of symbolic computation of the Hessian matrix. Symbolic computation often allows simplifications of the resulting expressions, which in turn means less conservative underestimators. We show by examples that even a small manipulation with the symbolic expressions, which can be processed automatically by computers, can have a large effect on the quality of underestimators.
\end{abstract}


\section{Introduction}

\subsection*{Convex underestimators}

To find a tight convex underestimator of an objective or/and constraint function is an essential problem in global optimization since it enables to easily compute a lower bound on the global optimal value, among others. In particular, it plays a crucial role in the well-known  global optimization $\alpha$BB method \cite{AdjDal1998, AdjAnd1998, AndMar1995, Flo2000, FloGou2009, FloPar2009}, which has been successfully applied in solving many real-life problems in biology and chemistry; see e.g. \cite{Flo2000,FloGou2009,SkjWes2012} and references therein. From the recent applications, let us mention utilizing of convex relaxations in biological systems \cite{MirPoz2012},  convexifications in semi-infinite programming \cite{ShiWu2012,Ste2012}, or application of convex relaxations in scheduling of crude oil operations \cite{LiMis2012}. See also the overview paper \cite{FloAkr2005}.

Let $f:\R^n\mapsto\R$ be a twice-differentiable objective or constraint function and $x_i\in\inum{x}_i=[\unum{x}_i,\onum{x}_i]$, $i=1,\dots,n$, interval domains for the variables. 
For symbolic manipulation we will also assume that we have an explicit analytic expression for $f(x)$, however, some basic improvement ideas from Section~\ref{sFurth} remain valid even for general case.
The aim is to construct a function $g:\R^n\mapsto\R$ satisfying:
\begin{enumerate}
\item
$f(x)\geq g(x)$ for every $x\in\ivr{x}$,
\item
$g(x)$ is convex on $x\in\ivr{x}$.
\end{enumerate}
The deterministic global optimization $\alpha$BB method \cite{AdjDal1998,AdjAnd1998, AndMar1995,Flo2000,FloPar2009} utilizes the convex underestimator in the form of
\begin{align}\label{fceUnder}
g(x):=f(x)-\sum_{i=1}^n\alpha_i(\onum{x}_i-x_i)(x_i-\unum{x}_i),
\end{align}
where $\alpha_i\geq0$, $i=1,\dots,n$, are determined such that $g(x)$ is convex. The Hessian of $g(x)$ reads
$$
\nabla^2 g(x)=\nabla^2 f(x)+2\diag(\alpha),
$$
where $\diag(\alpha)$ is the diagonal matrix with entries $\alpha_1,\dots,\alpha_n$. Let $\imace{H}$ be an interval matrix enclosing the image of $H(x):=\nabla^2 f(x)$ over $x\in\ivr{x}$. That is, the $(i,j)$th element of $\imace{H}$ is an interval $\inum{h}_{ij}=[\unum{h}_{ij},\onum{h}_{ij}]$ such that
$$
h_{ij}(x):=
\frac{\partial^2}{\partial x_i \partial x_j}f(x)\in\inum{h}_{ij},\quad
\forall x\in\ivr{x}.
$$
Now, to achieve convexity of $g(x)$, it is sufficient to choose $\alpha$ such that each matrix in $\imace{H}+2\diag(\alpha)$ is positive semidefinite, i.e., its eigenvalues are non-negative.
Eigenvalues of interval matrices were investigated e.g.\ in \cite{AdjDal1998,Flo2000,Hla2013a,HlaDan2011b,HlaDan2010,Mon2011}. For the purpose of the $\alpha$BB method, it seems that the most convenient method for bounding eigenvalues of interval matrices is the scaled Gerschgorin inclusion \cite{AdjDal1998,AdjAnd1998,Flo2000}. Its benefits are that it is easy to compute and eliminate the unknowns $\alpha_i$, $i=1,\dots,n$, and it is also usually sufficiently tight. For any positive $d\in\R^n$, we can put
\begin{align}\label{alpha}
\textstyle
\alpha_i:=\max\left\{0, -\frac{1}{2}
 \left(\unum{h}_{ii}-\sum_{j\not=i}|\inum{h}_{ij}|d_j/d_i\right)\right\},
\quad i=1,\dots,n,
\end{align}
where $|\inum{h}_{ij}|=\max\left\{|\unum{h}_{ij}|,|\onum{h}_{ij}|\right\}$. To reflect the range of the variable domains, it is recommended to use $d:=\onum{x}-\unum{x}$. Optimal choice of $d$ is discussed in \cite{Hla2013c}.

This classical $\alpha$BB approach was generalized in several ways. In \cite{AkrMey2004,SkjWes2012}, the authors considered convex underestimators in the form of
$$
g(x):=f(x)-(\onum{x}-x)^TP(x-\unum{x})+q,
$$
where $P\in\R^{n\times n}$ is a symmetric matrix with non-negative diagonal and $q\in\R$ is a correction value calculated so that the underestimation property is true. When $P$ is a diagonal matrix and $q=0$, the underestimator reduces to \nref{fceUnder}. 

Another class of underestimators defined as
$$
g(x):=f(x)-\sum_{i=1}^n 
 (1-e^{\gamma_i(\onum{x}_i-x_i)})(1-e^{\gamma_i(x_i-\unum{x}_i)})
$$
was discussed in \cite{AkrFlo2004,AkrFlo2004b,FloPar2009}, yielding the so called $\gamma$BB method. Herein, the parameters $\gamma_1,\dots,\gamma_n$ are computed via the classical $\alpha$BB method. Theoretical justification for  $\alpha$BB and $\gamma$BB relaxation terms is given in \cite{FloKre2007}.

Convex relaxations of quadratic functions were investigated in \cite{Ans2012}, linear relaxations in  \cite{DomNeu2012}, and a generalization of McCormick relaxations in \cite{ScoStu2011}. Another global optimization method, QBB, based on convex underestimators and branch \& bound scheme on simplices, was proposed in \cite{ZhuKun2005}.

\subsection*{Interval computation}

Interval computation \cite{HanWal2004,MooKea2009,Neu1990} serves to obtain rigorous enclosures to the image of intervals under various functions.
Let us introduce some notation.
An interval matrix $\imace{A}$ is defined as
$$
\imace{A}:=[\umace{A},\omace{A}]
=\{A\in\R^{m\times n}\mmid \umace{A}\leq A\leq\omace{A}\},
$$
where $\umace{A},\omace{A}\in\R^{m\times n}$ are given. 
The center and radius of $\imace{A}$ are respectively defined as
$$
\Mid{A}:=\frac{1}{2}(\umace{A}+\omace{A}),\quad
\Rad{A}:=\frac{1}{2}(\omace{A}-\umace{A}).
$$
Interval vectors and intervals can be regarded as special interval matrices of sizes $m$-by-$1$ and $1$-by-$1$, respectively. 

Let $f:\R^n\mapsto \R$ and an interval vector $\ivr{x}$ be given. The image
$$
f(\ivr{x}):=\{f(x)\mmid x\in\ivr{x}\}
$$
is hard to determine in general. That is why one usually seeks for its enclosure, i.e., an interval $\inum{f}$ such that $f(\ivr{x})\subseteq\inum{f}$.
Interval arithmetic extends the standard arithmetic naturally as follows. 
Let $\inum{a}=[\unum{a},\onum{a}]$ and $\inum{b}=[\unum{b},\onum{b}]$ be intervals, then we define
\begin{align*}
\inum{a}+\inum{b}&=[\unum{a}+\unum{b},\onum{a}+\onum{b}],\\
\inum{a}-\inum{b}&=[\unum{a}-\onum{b},\onum{a}-\unum{b}],\\
\inum{a}\inum{b}&=
 [\min(\unum{a}\unum{b},\unum{a}\onum{b},\onum{a}\unum{b},\onum{a}\onum{b}),
  \max(\unum{a}\unum{b},\unum{a}\onum{b},\onum{a}\unum{b},\onum{a}\onum{b})],\\
\inum{a}/\inum{b}&=
 [\min(\unum{a}/\unum{b},\unum{a}/\onum{b},\onum{a}/\unum{b},\onum{a}/\onum{b}),
  \max(\unum{a}/\unum{b},\unum{a}/\onum{b},\onum{a}/\unum{b},\onum{a}/\onum{b})],
\end{align*}
with $0\not\in\inum{b}$ in case of division. The image of an interval for the basic functions such as sine, cosine, exponential can be determined by a direct inspection. Thus, by using interval arithmetic, we can evaluate many algebraic expressions on intervals. However, notice two points. First, the results may be highly overestimated, and, second, two mathematically equivalent expressions may yield enclosures of different quality.

For example, consider a trivial example
$$
f=(x-3)^2=x^2-6x+9.
$$
and $x\in\inum{x}=[1,4]$.
Evaluating $(\inum{x}-3)^2$ gives $[0,4]$, but $\inum{x}^2-6\inum{x}+9=[-14,19]$. Therefore symbolical manipulation of expressions in order to make then as simple as possible may dramatically influence tightness of the calculated enclosure. This principle is highlighted in this paper, and confirmed by examples.

Besides interval arithmetic, there are other methods to compute enclosures of the function images on intervals. For instance, by utilizing the mean value theorem, we obtain the so called mean value form of function enclosure. For simplicity, let $f:\R\mapsto\R$ be univariate, $\inum{x}$ an interval and $a\in\inum{a}$. Then
$$
f(\ivr{x})\subseteq f(a)+f'(\ivr{x})(\inum{x}-a),
$$
where $f'(\ivr{x})$ is an enclosure to the derivative of $f$ on $\inum{x}$. For a generalization to multivariate case see e.g.\ \cite{HanWal2004,MooKea2009,Neu1990}. The performance of mean value form can be improved by replacing derivatives by slopes. The slope of $f$ at $a\in\inum{x}$ is defined as
$$
S_f(x,a):=
\begin{cases}
\frac{f(x)-f(a)}{x-a} &\mbox{if }  x\not=a,\\
f'(x) & \mbox{otherwise}.
\end{cases}
$$
Slopes can be evaluated in a similar manner as derivatives, but the result provably outperforms derivatives. Moreover, slopes can handle also some non-smooth functions such as the absolute value (which is convenient in our approach). For more details, see e.g.\ \cite{HanWal2004,MooKea2009,Neu1990}.

\section{Symbolic computation of $\alpha$}

In this section, we study computation of $\alpha$ from \nref{alpha} and its impact on the quality of convex underestimators for the classical $\alpha$BB  method. 


The proposed idea behind more effective computation of $\alpha$ is to directly substitute for the Hessian entries in that formula instead of computing an interval enclosure of the Hessian and then using those entries.

Define
\begin{align}\nonumber
h_i(x)&:=
\frac{\partial^2}{\partial x_i^2}f(x)-\sum_{j\not=i}
  \left|\frac{\partial^2}{\partial x_i \partial x_j}f(x)\right|d_j/d_i\\
&=h_{ij}(x)-\sum_{j\not=i}\left|h_{ij}(x)\right|d_j/d_i,
\quad i=1,\dots,n.\label{hi}
\end{align}
The entries of $\alpha$ then follows
\begin{align}\label{alphaHi}
\alpha_i:=\max\left\{0, -\frac{1}{2}\,\unum{h_i(\ivr{x})} \right\},
\quad i=1,\dots,n,
\end{align}
If we compute the images $h_i(\inum{x})$ by using interval arithmetic and automatic differentiation, the result will be the same as for the classical case. However, if we employ symbolic differentiation and rearrangements of the expressions, the overall overestimation can considerably be reduced.

First notice that provided $\inum{h}_{ij}$ does not include zero in its interior, then the sign of
$h_{ij}(x)$
is stable (invariant) and we can remove the corresponding absolute value in \nref{hi}. Provided $h_{ij}(x)$ are sign stable for all $j\not=i$, the function $h_i(x)$ is found continuous (and differentiable if $f(x)$ is higher order differentiable), and thus tighter enclosure of the image $h_i(\inum{x})$ can be expected by using appropriate interval methods (monotonicity checking \cite{Han1997} etc.). In principle, even when some of the terms $h_{ij}(x)$ are recognized as sign stable, we may achieve good results.

Now suppose that the Hessian matrix $\nabla^2 f(x)$ is computed symbolically. Thus, we have an explicit formula for $h_i(x)$ and an enclosure of its image can be calculated not only by interval arithmetic, but also by the mean value form using slopes or any other suitable technique. Moreover, we can symbolically manipulate and rearrange the formula for $h_i(x)$ in order to achieve a more convenient form for interval evaluation. In the next section, we demonstrate by several examples that even a simple expression rearrangement, which can be done automatically by computer, may result in large increase of performance.

\section{Computational studies}

We present some numerical experiments done in \textsf{MATLAB}, and we employed the interval toolbox \textsf{INTLAB v6} \cite{Rum1999}. The toolbox provides us with the interval arithmetic, images of basic functions over intervals, interval gradients and interval Hessian matrices. Notice that in the examples below, the vector $\alpha$ computed by the $\alpha$BB method may slightly differ from the literature values just because we calculated the initial interval Hessian numerically by \textsf{INTLAB}.

\begin{example}
Consider the function from \cite{GouFlo2008,SkjWes2012}
\begin{align*}
f(x_1,x_2,x_3,x_4)
=(x_1+10x_2)^2+5(x_3-x_4)^2+(x_2-2x_3)^4+10(x_1-x_4)^4,
\end{align*}
where $x\in\ivr{x}=[0,1]^4$. It is known that the global minimum is $f^*=0$.

First, we compute the interval Hessian
$$
\nabla^2 f(\ivr{x})\subseteq \imace{H}=
\begin{pmatrix}
[ -118,  122]& [   20,   20]& [    0,    0]& [ -120,  120] \\{}
[   20,   20]& [  176,  248]& [  -96,   48]& [    0,    0]\\{}
[    0,    0]& [  -96,   48]& [  -86,  202]& [  -10,  -10] \\{}
[ -120,  120]& [    0,    0]& [  -10,  -10]& [ -110,  130] 
\end{pmatrix}.
$$
By the scaled Gerschgorin method we obtain
$$
\alpha=(129,0,96,120).
$$
and the corresponding lower bound on $f^*$ is $-85.1312$.

Let us compute the Hessian matrix symbolically
$$
\nabla^2 f({x})=
\begin{pmatrix}
2+120(x_1-x_4)^2 & 20 & 0 & -120(x_1-x_4)^2 \\
20 & 200+12(x_2-2x_3)^2 & -24(x_2-2x_3)^2 & 0\\
0 & -24(x_2-2x_3)^2 & 10+48(x_2-2x_3)^2 & -10\\
-120(x_1-x_4)^2 & 0 & -10 & 10+120(x_1-x_4)^2
\end{pmatrix}.
$$
Since all off-diagonal entries are sign stable, we can omit the absolute values in \nref{hi}. The function
$$
h_1(x)=2+120(x_1-x_4)^2 -20-120(x_1-x_4)^2
$$
is evaluated by interval arithmetic with the result $[-138,102]$, so we put $\alpha_1=69$. Analogously we proceed further and get
$$
\alpha=(69,0,48,60).
$$
The corresponding lower bound on $f^*$ is $-43.2171$.

However, we can obtain yet much tighter lower underestimator. Simplifying $h_1(x)$ to $h_1(x)=-18$, and similarly for the others, we have
$$
\alpha=(18,0,0,0)
$$
and the lower bound on the global minimum is $-1.9768$.
\end{example}

\begin{figure}[t]
\begin{minipage}[b]{.44\linewidth}
\begin{center}
\pgfimage[width=1.0\linewidth]{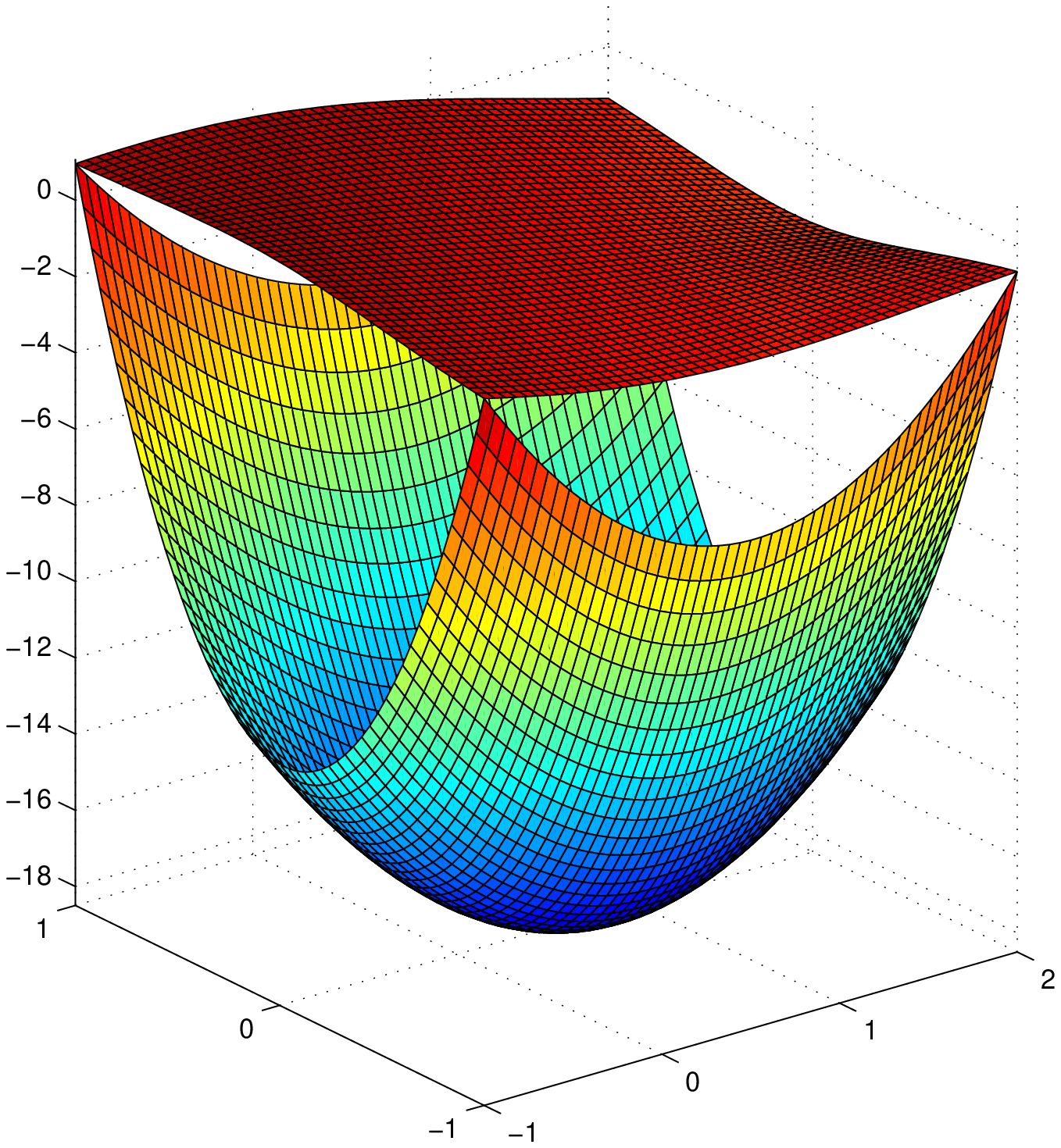}
\caption{(Example~\ref{ex1}) The function and its convex underestimator for the classical $\alpha$BB method.\label{fig-ex-1a}}
\end{center}
\end{minipage}
\hfill 
\begin{minipage}[b]{.44\linewidth}
\begin{center}
\pgfimage[width=1.0\linewidth]{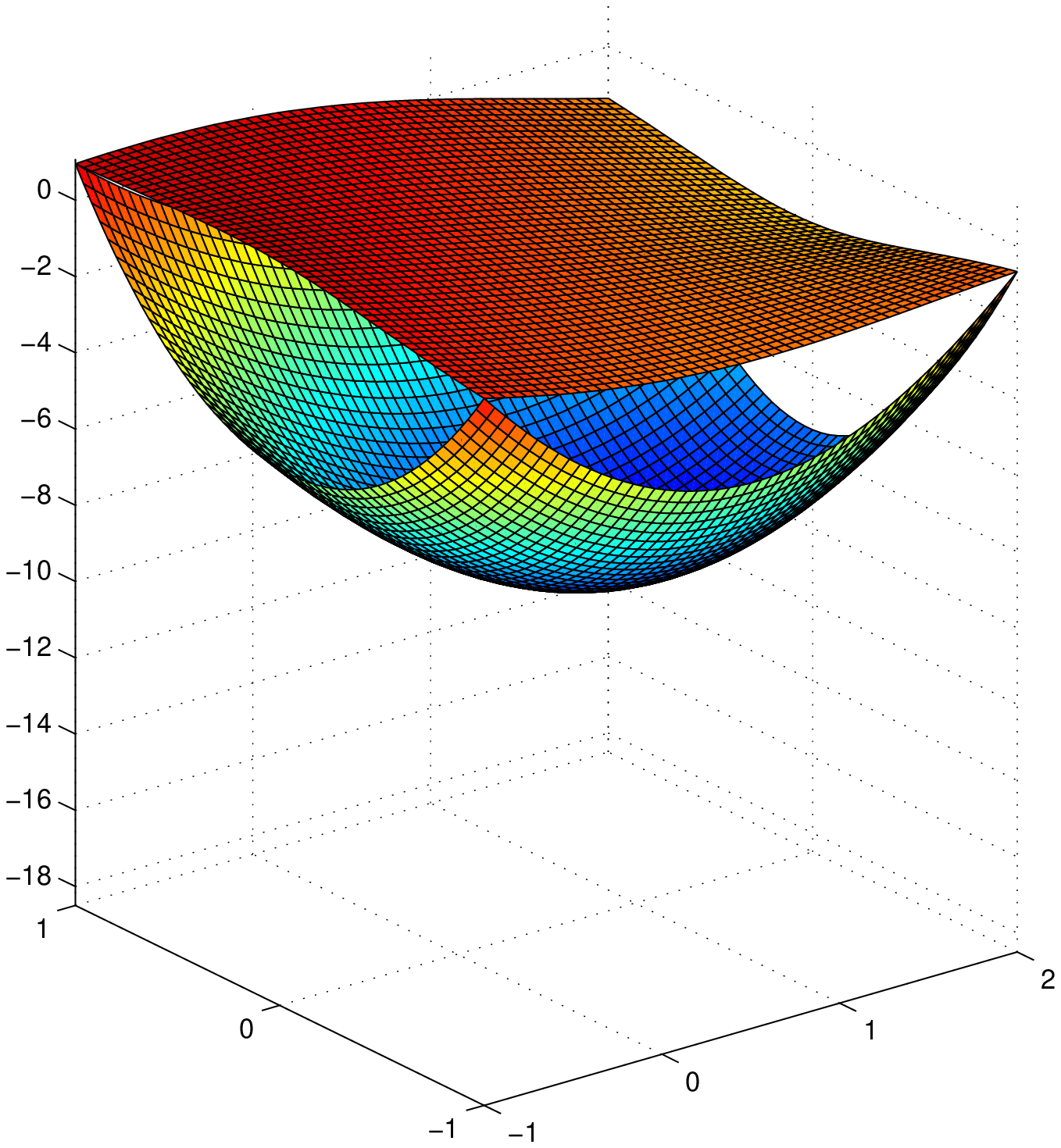}
\caption{(Example~\ref{ex1}) The function and its convex underestimator by our approach.\label{fig-ex-1b}}
\end{center}
\end{minipage}
\end{figure}

\begin{example}\label{ex1}
Consider the function from \cite{AdjDal1998,Flo2000}
\begin{align*}
f(x_1,x_2)
=\cos(x_1)\sin(x_2)-\frac{x_1}{x_2^2+1},
\end{align*}
where $x_1\in[-1,2]$ and $x_2\in[-1,1]$. The optimal value is known to be $f^*=-2.02181$.

Proceeding along the classical $\alpha$BB method, we compute 
$$
\nabla^2 f(\ivr{x})\subseteq \imace{H}=
\begin{pmatrix}
[   -0.8415,  0.8415]& [   -5.0000,    4.8415] \\{}
[   -5.0000,  4.8415]& [  -18.8415,   20.8415]  
\end{pmatrix},
$$
whence 
$$
\alpha=(2.0874,13.1707),
$$
and the corresponding lower bound on $f^*$ is $-18.4970$.

Using the symbolical approach, we express the Hessian matrix as in \cite{AdjDal1998,Flo2000}
$$
\nabla^2 f({x})=
\begin{pmatrix}
-\cos(x_1)\sin(x_2) &
\displaystyle-\sin(x_1)\cos(x_2)+\frac{2x_2}{(x_2^2+1)^{2}} \\[1em]
\displaystyle-\sin(x_1)\cos(x_2)+\frac{2x_2}{(x_2^2+1)^{2}} &
\displaystyle-\cos(x_1)\sin(x_2)+
 \frac{2x_1(x_2^2+1)^2-8x_1^{}x_2^2(x_2^2+1)}{(x_2^2+1)^{4}}
\end{pmatrix},
$$
and have to evaluate the functions
\begin{align*}
h_1(x)&=-\cos(x_1)\sin(x_2)
 -\frac{2}{3}\left|-\sin(x_1)\cos(x_2)+\frac{2x_2}{(x_2^2+1)^{2}}\right|,\\
h_2(x)&=-\cos(x_1)\sin(x_2)+
 \frac{2x_1(x_2^2+1)^2-8x_1^{}x_2^2(x_2^2+1)}{(x_2^2+1)^{4}}
 -\frac{2}{3}\left|-\sin(x_1)\cos(x_2)+\frac{2x_2}{(x_2^2+1)^{2}}\right|.
\end{align*}
We cannot get rid of the absolute values since the off-diagonal entries of the Hessian are not sign stable.
The direct evaluation of the function thus makes no improvement, but we can easily simplify the expression for $h_2(x)$,
\begin{align*}
h_2(x)&=-\cos(x_1)\sin(x_2)+
 \frac{2x_1(2-6x_2^2)}{(x_2^2+1)^{3}}
 -\frac{2}{3}\left|-\sin(x_1)\cos(x_2)+\frac{2x_2}{(x_2^2+1)^{2}}\right|.
\end{align*}
Now, we calculate
$$
\alpha=(1.4208,5.4208),
$$
and the lower bound on the optimal value is $-9.3110$.
\end{example}

\begin{figure}[t]
\begin{minipage}[b]{.44\linewidth}
\begin{center}
\pgfimage[width=1.0\linewidth]{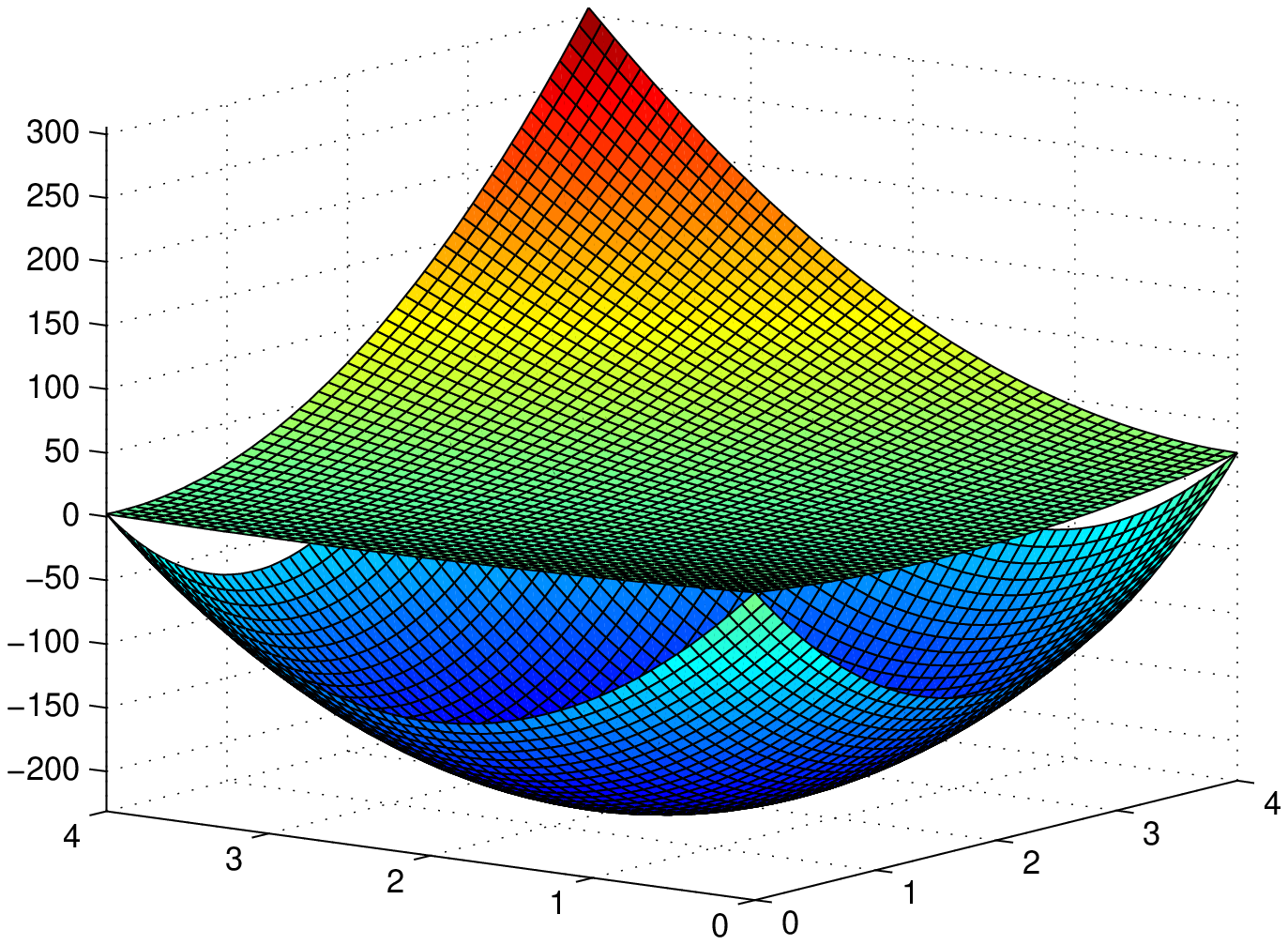}
\caption{(Example~\ref{ex8}) The function and its convex underestimator for the classical $\alpha$BB method.\label{fig-ex-8a}}
\end{center}
\end{minipage}
\hfill 
\begin{minipage}[b]{.44\linewidth}
\begin{center}
\pgfimage[width=1.0\linewidth]{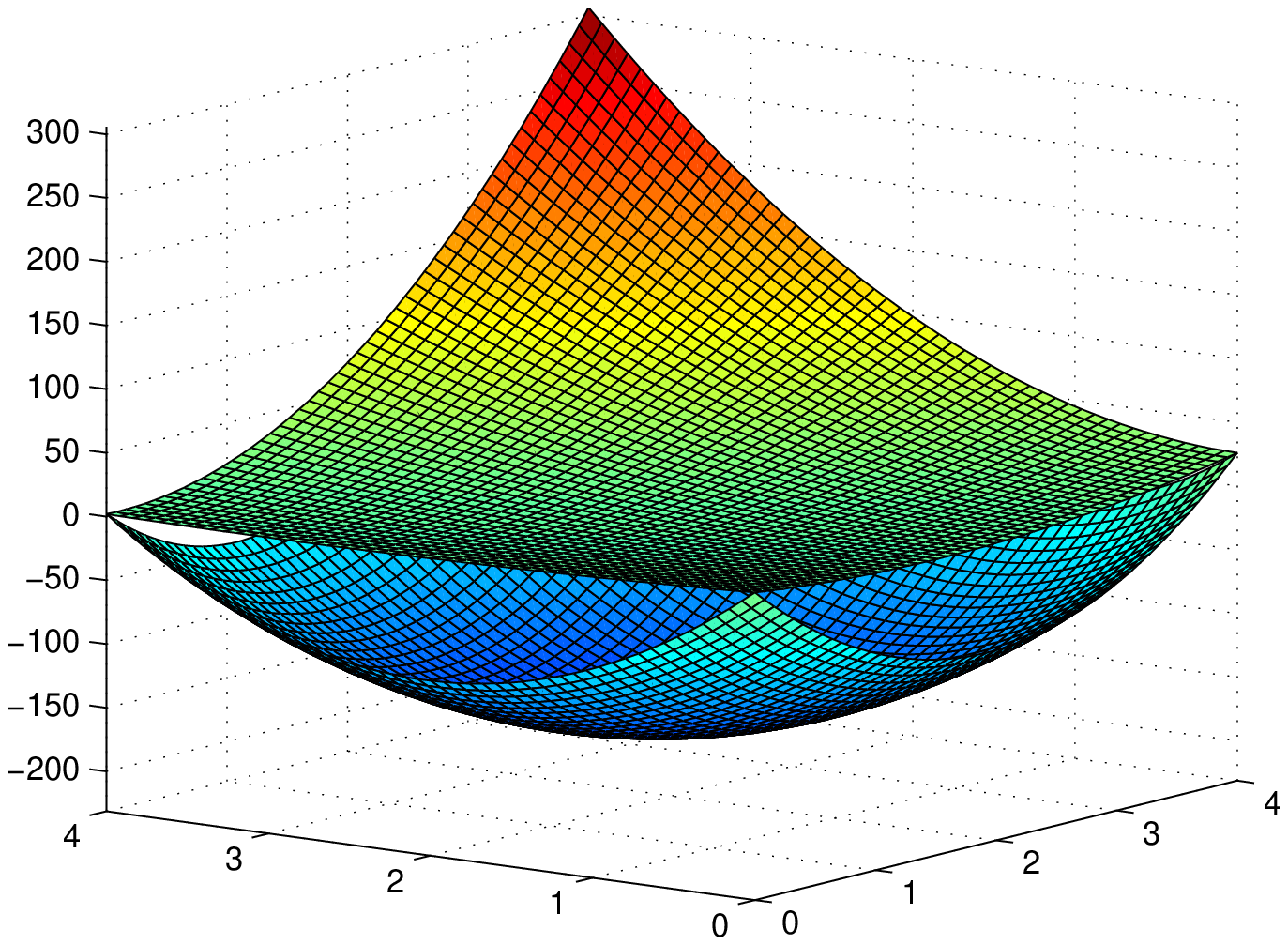}
\caption{(Example~\ref{ex8}) The function and its convex underestimator by our approach.\label{fig-ex-8b}}
\end{center}
\end{minipage}
\end{figure}

\begin{example}\label{ex8}
Consider the function from \cite{SkjWes2012}
\begin{align*}
f(x_1,x_2)
=(2x_1+x_2-3)^2+(x_1x_2-1)^2,
\end{align*}
where $x\in[0,4]^2$. The optimal value is $f^*=0$. 

The classical $\alpha$BB method computes
$$
\alpha=(29,32),
$$
and the lower bound on $f^*$ is  $-231.0459$.
The generalization of the $\alpha$BB method using non-diagonal quadratic terms improves the lower bound only to $-230.90$.

Evaluating the Hessian matrix symbolically and the functions $h_1(x)$ and $h_2(x)$ by the mean value form, we obtain
$$
\alpha=(21,24),
$$
and arrive at the lower bound  $ -168.1901$. Thus, we tighten the lower bound by $27.2\%$ without using any algebraic simplifications of the Hessian or the functions $h_i(x)$.
\end{example}

\begin{figure}[t]
\begin{minipage}[b]{.44\linewidth}
\begin{center}
\pgfimage[width=1.0\linewidth]{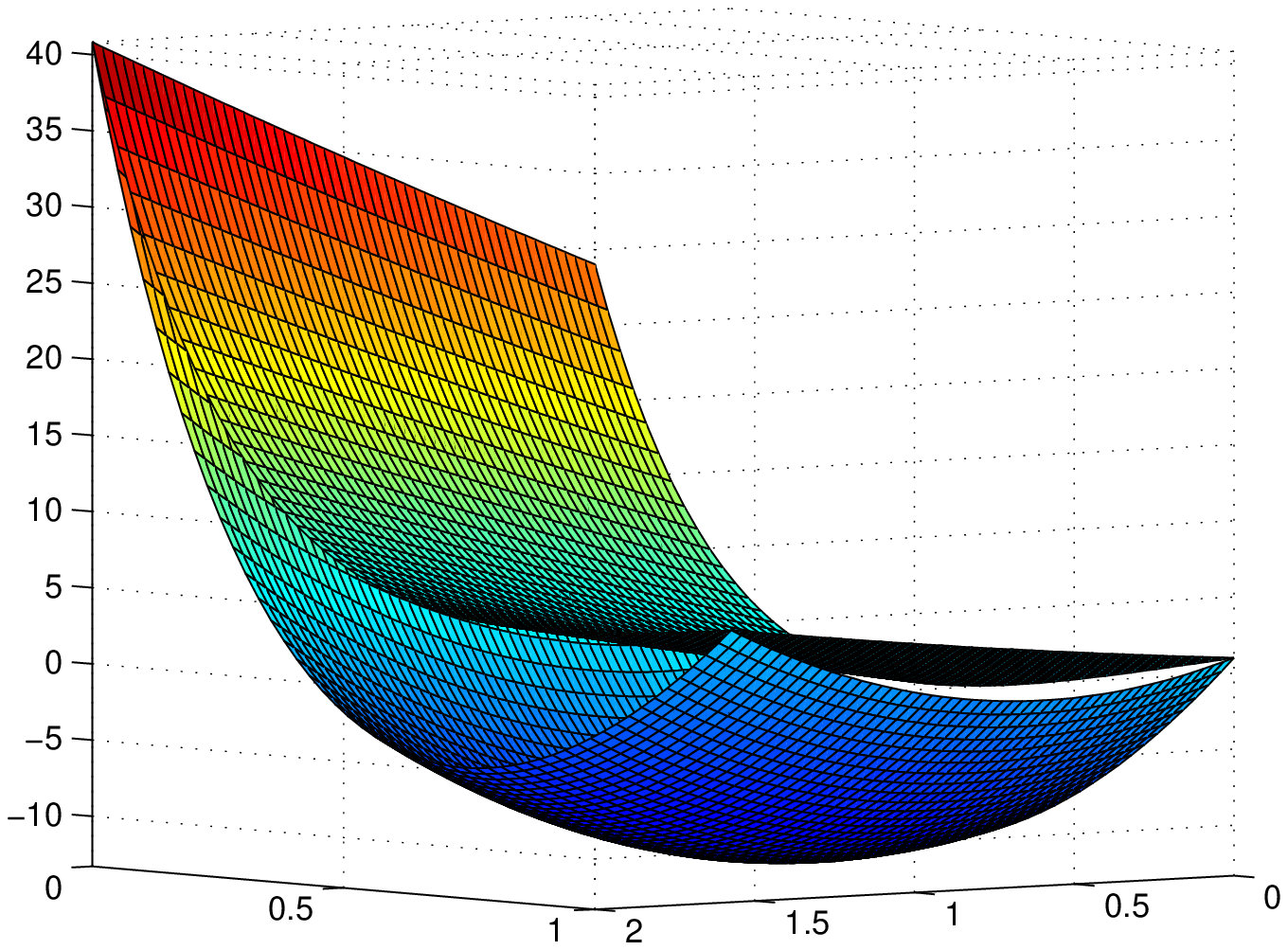}
\caption{(Example~\ref{ex9}) The function and its convex underestimator for the classical $\alpha$BB method.\label{fig-ex-9a}}
\end{center}
\end{minipage}
\hfill 
\begin{minipage}[b]{.44\linewidth}
\begin{center}
\pgfimage[width=1.0\linewidth]{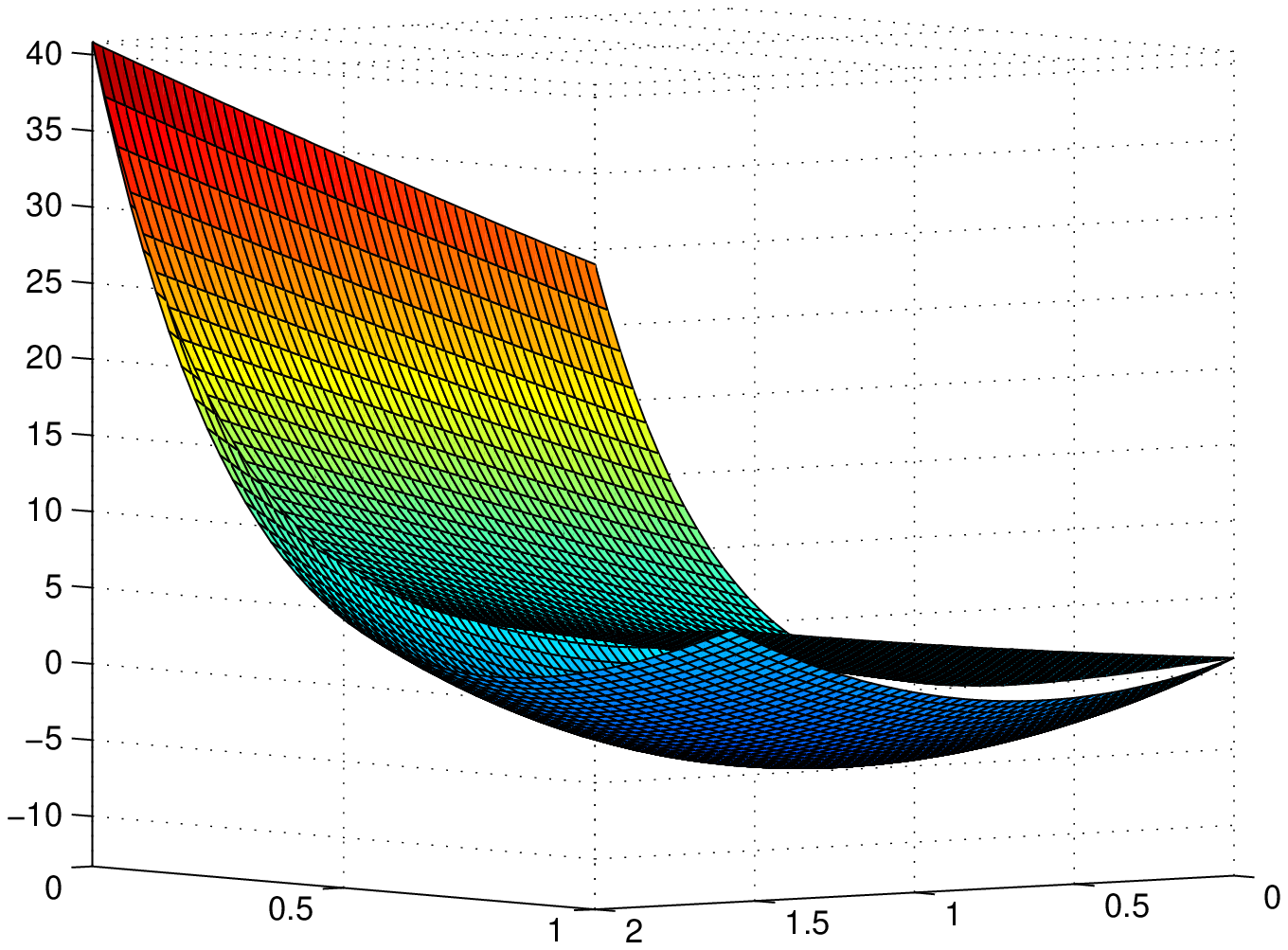}
\caption{(Example~\ref{ex9}) The function and its convex underestimator by our approach (after simplification of $h_2(x)$).\label{fig-ex-9b}}
\end{center}
\end{minipage}
\end{figure}

\begin{example}\label{ex9}
Consider the function from \cite{SkjWes2012}
\begin{align*}
f(x_1,x_2)
=(1+x_1-e^{x_2})^2,
\end{align*}
where $x_1\in[0,1]$ and $x_2\in[0,2]$. The optimal value is $f^*=0$. 

According to \cite{SkjWes2012}, the classical interval  $\alpha$BB method computes the lower bound on the optimal value as $-14.46$, which was improved by the non-diagonal generalization to $-13.18$. In  \cite{SkjWes2012}, the authors also determined the optimal enclosure to the interval Hessian, which resulted in the lower bound $-8.24$ by the classical method and to $-6.94$ by the non-diagonal generalization.

In our approach, we handle the functions
\begin{align*}
h_1(x)&=2-4e^{x_2},\\
h_2(x)&=2e^{x_2}e^{x_2}-2(1+x_1-e^{x_2})e^{x_2}-e^{x_2}.
\end{align*}
By using interval arithmetic or mean values form, we obtain the lower bound $-12.65$. Evaluation of $h_1(x)$ on intervals is always optimal since each variable occurs at most once in the expression. The expression for $h_2(x)$ is easy to simplify to 
\begin{align*}
h_2(x)&=(-3-2x_1+4e^{x_2})e^{x_2}.
\end{align*}
Now, the resulting lower bound is $-6.5629$.
\end{example}

\section{Further improvements}\label{sFurth}

Recall that as long as  $\inum{h}_{ij}$ does not include zero in its interior, then the sign of $h_{ij}(x)$ is stable and we can remove the corresponding absolute value in \nref{hi}. The computational studies presented in the previous section showed that this enables us to compute the image of $h_i(x)$ more efficiently. 
When $\inum{h}_{ij}$ includes zero in its interior, the evaluation of $h_i(\ivr{x})$ is more challenging. Let us discuss some more promising approaches than the direct evaluation by interval arithmetic is.

We do not need to determine a tight enclosure to the whole image $h_i(\ivr{x})$, but in view of \nref{alphaHi} only a tight lower bound on $\unum{h_i(\ivr{x})}$. This means that we can estimate $|h_{ij}(x)|$ from above.

Assume without loss of generality that $\unum{h}_{ij}+\onum{h}_{ij}\geq0$, otherwise we consider $-h_{ij}(x)$ instead of $h_{ij}(x)$. Then  $|h_{ij}(x)|\leq h_{ij}(x)-\unum{h}_{ij}$ disposes the absolute value. Using this estimation may or may not result in a tighter enclosure. However, provided $\unum{h}_{ij}$ is close to the zero, we can expect that this estimation is effective, or at least the worsening is very small (always bounded by $\unum{h}_{ij}$).

Another possibility is to estimate the absolute value from above by the tightest linear function \cite{Bea1998}. 

\begin{proposition}\label{thmBea}
For every $y\in\inum{y}\subset\R$ with $\unum{y}<\onum{y}$ one has
\begin{align}\label{thmBeaAbs}
|y|\leq \gamma y+\beta,
\end{align}
where
\begin{align*}
\gamma=\frac{|\onum{y}|-|\unum{y}|}{\onum{y}-\unum{y}}\ \mbox{ and }\ 
\beta=\frac{\onum{y}|\unum{y}|-\unum{y}|\onum{y}|}{\onum{y}-\unum{y}}.
\end{align*}
Moreover, if $\unum{y}\geq0$ or $\onum{y}\leq0$ then \nref{thmBeaAbs} holds as equation.
\end{proposition}

Employing this proposition, we simply estimate
\begin{align}\label{ineqBeaAbsH}
|h_{ij}(x)|\leq  \gamma h_{ij}(x)+\beta,
\end{align}
where
\begin{align*}
\gamma=\frac{|\onum{h}_{ij}|-|\unum{h}_{ij}|}{\onum{h}_{ij}-\unum{h}_{ij}}\ \mbox{ and }\ 
\beta=\frac{\onum{h}_{ij}|\unum{h}_{ij}|-\unum{h}_{ij}|\onum{h}_{ij}|}{\onum{y}-\unum{h}_{ij}}.
\end{align*}
Since \nref{thmBeaAbs} is the best linear upper approximation of the absolute value, this relaxation can never be worse than the direct interval evaluation of $|h_{ij}(\ivr{x})|$, since it estimates the value of the function by the constant $\onum{h}_{ij}$. In contrast, linear relaxation of $|h_{ij}(\ivr{x})|$ by means of \nref{ineqBeaAbsH} is suitable for symbolic simplifications of $h_i(x)$.

\begin{example}
Consider the function 
\begin{align*}
f(x_1,x_2)
=20x_1x_2^2+10x_1^3-4x_3^3-7x_1^2-70x_1x_2
\end{align*}
where $x_1,x_2\in[1,2]$. Its Hessian matrix reads
\begin{align*}
\nabla^2 f(x)=
\begin{pmatrix}
60x_1-7 & 40x_2-70 \\  40x_2-70 & 40x_1-24x_2
\end{pmatrix}.
\end{align*}
Evaluation by interval arithmetic leads to the interval enclosure
\begin{align*}
\nabla^2 f(\ivr{x})\subseteq\imace{H}=
\begin{pmatrix}
[53, 113] & [-30,10] \\  [-30,10] &  [-8,56]
\end{pmatrix}.
\end{align*}
The classical computation of $\alpha$ by \nref{alpha} results in 
$
\alpha=(0,19).
$
Let us compare it with the proposed two ways to relax the absolute value. First, we estimate
$$
|h_{21}(x)|
 = |70-40x_2|
 \leq  h_{21}(x)-\unum{h}_{21}
 = 70-40x_2 +10=80-40x_2.
$$
Now, we calculate 
$$
h_2(x)\geq 40x_1-24x_2 - (80-40x_2)=40x_1+16x_2 -80\in[-24,32],
$$
whence $\unum{h_2(\ivr{x})}\geq -24$, and therefore $\alpha_2=12$ is notably tightened.

In the second way, we compute the coefficients $\gamma=-0.5$ and $\beta=15$ corresponding to $\inum{h}_{21}$. This leads to the estimation
\begin{align*}
|h_{21}(x)|=| 40x_2-70 |\leq  -0.5(40x_2-70)+15=-20x_2+50.
\end{align*}
Thus,
$$
h_2(x)\geq 40x_1-24x_2 -(-20x_2+50)=40x_1-4x_2 -50\in[-18,26],
$$
and we get yet lower value of $\alpha_2=9$.
\end{example}

\section{Conclusion}

We presented a variant of the convex underestimator construction in the $\alpha$BB method. We discussed the advantages of computing the Hessian matrix symbolically. Compared to automatic differentiation, we can utilize various techniques from interval computation area to obtain tighter results. The numerical experiments demonstrated that only a small symbolic simplification of expressions may have a large effect on the quality of the resulting underestimators. 

A function maybe expressed by using many equivalent algebraic formulae. It is not always clear which one to choose for interval evaluation. However, as shown by our examples, even a small rearrangement can yield much tighter underestimators than other generalizations and improvements of the $\alpha$BB method. Therefore, we recommend to pay more attention to symbolic handling with expressions and drive the research in this direction.

\subsubsection*{Acknowledgments.} 

The author was supported by the Czech Science Foundation Grant P402/13-10660S.


\bibliographystyle{abbrv}
\bibliography{conv_under}

\end{document}